\documentclass[9pt]{article}
\usepackage{pdflscape}
\usepackage{eurosym}
\usepackage{soul}
\usepackage{calc}
\usepackage{color}
\usepackage{amsfonts}
\usepackage{latexsym}
\usepackage{placeins}
  \usepackage[dvips]{graphicx}
\usepackage{amssymb}
\usepackage{bbold}
\usepackage{authblk}
\usepackage{amsmath}
\usepackage[utf8]{inputenc}
\usepackage{eurosym}
\usepackage{multirow}

\usepackage{setspace}

\newcommand{\smalllineskip}{\baselineskip=15pt}

\renewenvironment{abstract}[0]{\small\rm
        \begin{center}ABSTRACT
        \\ \vspace{1pt}
        \begin{minipage}{5.2in}\smalllineskip
        \hspace{1pc}}{\end{minipage}\end{center}\vspace{-1pt}}

\title{Plant Performance in Precision Horticulture:\\Optimal climate control under stochastic dynamics }
\date{}

\author[1]{Simon van Mourik} \author[2]{Bert van 't Ooster}  \author[3]{Michel Vellekoop}
\affil[1,2]{Plant Sciences Group,
Wageningen University, The Netherlands}
\affil[3]{Faculty of Economics \& Business,
University of Amsterdam, The Netherlands}

\begin{document}
\maketitle
\centerline{\today}

\newcommand{\mycomment}[1]{}
\newcommand{\be}{\begin{eqnarray*}}
\newcommand{\ee}{\end{eqnarray*}}
\newcommand{\mycomments}[1]{}
\newcommand{\quav}[2]{\left< #1,#2 \right>}
\newcommand{\sfrac}[2]{ {\textstyle\frac{#1}{#2}} }

\bibliographystyle{unsrt}

\begin{abstract}


This paper presents a risk mitigating, time-varying feedback control algorithm for crop production when state dynamics are subject to uncertainty. The model based case study concerns a 40 day production round of lettuce in a greenhouse where control input consists of daily and nightly temperature set points. The control problem was formulated in terms of a stochastic Markov decision process with the objective to maximize the expected net revenue at harvest time.
 
The importance of time-varying feedback and of risk mitigation was investigated by making a comparison with a  controller that takes uncertainty into account but is static and a controller which is dynamic but ignores the uncertainty in the state dynamics. For the case of heat limited crop growth, and strict requirements on harvest weight precision, the dynamic stochastic controller outperformed the static controller in terms of both maximal expected net revenue (by 19 \%) and state precision at harvest time (with 50 \% less standard deviation). It also outperformed the deterministic controller for both criteria (15 \% in maximal expected net revenue and 8 \% less standard deviation). 
A detailed sensitivity analysis showed that such improvements in performance levels are robust, since they hold over large ranges of uncertainty in state dynamics, required harvest precision levels, starting days, and initial weights. 

The results provide insights in potential of dynamic feedback and risk mitigation strategies for high precision requirements under uncertainty. Although the results should be interpreted with caution, they illustrate the considerable potential benefit for stochastic greenhouse climate control under uncertainty when high precision is required.

\end{abstract}

\section*{Keywords}
Markov decision process, greenhouse, lettuce
 
\section{Introduction}

In precision horticulture, there are various causes that lead to uncertainty in the net revenues. First off, there exists quite some unexplained variation in crop development, which is corroborated by the observation that crop models display a broad range in yield prediction accuracy, from $1\%$ to $15\%$ \cite{Poorter2013,Righini2020}. Uncertainty may also increase because of spatiotemporal variability in climatic variables \cite{balendonck2010} that is unforeseen or not modelled. Second, due to sensor noise it might be hard to obtain accurate estimates of the state of the crop and the climate, that form the starting point of any predictions about future development. Filter algorithms have been shown to be able to mitigate the effect of sensing errors using model based predictions \cite{Speetjens2009,Hameed2010}. However, even well tuned filters may not meet the accuracy standards that are required in modern greenhouses \cite{vanMourik2019}. Third, forecast errors of processes like the weather, energy prices and crop market can substantially increase prediction uncertainty and thereby hinder control performance \cite{Kuijpers2021,Su2021,Payne2021}. Uncertainty in crop response and indoor climate dynamics is often modelled with parametric uncertainty \cite{Vazquez2014,Lopez2018uncertainty, Guzman2009}, or stochastic noise \cite{Speetjens2009}. State estimation methods such as Kalman filters are based on modeling the stochastic processes for noise in the state and the measured output \cite{Kuijpers2021,Oldewurtel2013}.  

In this paper we focus on the control problem of optimizing the expected costs of a dynamic system with stochastic uncertainty. Our problem is in the class of Markov decision processes (MDP) \cite{boucherie2017} and can be solved by stochastic dynamic programming techniques \cite{Ross2014}. Research that concerns MDP in precision farming spans a wide range of applications such as precision field irrigation \cite{Huong2018}, pig breeding \cite{kovacs2012}, crop disease control \cite{Onstad1985,sells1995}, and planting, fertilization and harvesting \cite{Boussios2019}. Greenhouse horticulture research for optimal control under stochastic disturbances includes optimal indoor climate management using a stochastic error forecast model \cite{Kuijpers2022, Della2019}, optimal energy management under random renewable energy availability \cite{Zhuang2018}, temperature control under stochastic energy prices \cite{garcia2023multi}, and tracking and predicting crop states such as dry weight and leaf area in lettuce using a dynamic Bayesian network \cite{Kocian2020}.

To deal with the high dimensional decision space and the corresponding computational challenges that characterize Markov decision process optimization, reinforcement learning for greenhouse climate control is an active topic of investigation \cite{Tchamitchian2005, zhang2021, Ajagekar2022}. Stochastic dynamics add a layer of complexity to the control problem, which increases computational demand and might hinder transparency of software and result in outcomes that are harder to explain. The question that arises is: what is the added performance of a complex stochastic controller compared to more straightforward control methods? So far, not much has been reported about the effectiveness of stochastic control in a greenhouse setting. A study which found that a stochastic controller outperformed expert growers in a tomato greenhouse regarding net profit (93\%) and crop yield (by 10\%) \cite{cao2022igrow}, suggests substantial benefits, but the question remains to what extent the complexity of the controller contributes to the performance. 

This paper introduces a dynamic stochastic controller that is obtained by optimizing an MDP in terms of expected net revenue by means of dynamic programming. The case study concerns a 40 day production round of lettuce in a greenhouse where control input consists of daily and nightly temperature set points. The control problem was formulated  with the objective to maximize the expected net revenue at harvest time and we obtain a control strategy with an excellent performance, a clear interpretation and limited complexity. 

To study the properties of the optimal control strategy, the control policy, and corresponding dynamics and performance are visualized in order to gain insight into the rationale behind the control policy. The results are visualized in four maps that show on the discrete grid of time and state the controller's set points for day and night temperatures, the associated crop growth rates, and the associated value functions in terms of expected net revenue.  

Subsequently, the optimized dynamic stochastic time-varying feedback controller is compared against a static stochastic controller and a dynamic but deterministic time-varying feedback controller. Performance was assessed in terms of maximum expected net revenue at the start of a production round, the crop state uncertainty at harvesting time, and robustness against deviations from the optimal starting time.

\section{Materials and methods}

\subsection{System dynamics}

For times $k\in {\mathbb K}:=\{0,1,2,3...,T-1\}$ we define the stochastic scalar dynamic system 
\begin{align}\label{eq:model}
x_{k+1}&= x_k+f_k(x_k,u_k)\,\Delta t  \ + \ x_k\sigma\sqrt{\Delta t} \ \epsilon_k\\
u_k &= g_k(x_k).
\end{align}
Here $x$ is the state variable, with time index $k$;  time is discretized as $t_k=k\Delta t$ for a given $\Delta t>0$. The function $f:{\mathbb R}^3\to{\mathbb R}$ describes the dynamic interactions between state and input. The vectors $u_k$ are the control inputs, which are determined by a state feedback control policy $g$, which consists of a sequence of input functions $g_k$ for all time steps $k$. The stochastic variables $(\epsilon_k)_{k\in \mathbb K}$ are independent and identically normally distributed. We use the notation $x^g_k$ for the states that are generated if a policy $g$ is applied for a given initial state $x_0^g=x_0$.

\subsection{Control problem}
The goal is to determine a policy $g$ that maximizes value, defined by the expected revenues  $J(x_T^g)$ at the final time $T$, minus the sum over all running costs $L_k$. The associated control problem is to find a policy that solves the optimization problem
\begin{equation}\label{eq:controlproblem_0}
\max_{g\in {\cal G}} \mathbb{E} \left[ J(x_T^g) - \sum_{k=0}^{T-1} L_k(g_k(x_k^g))) \Delta t\, \right]
\end{equation} 
 over all  control laws $(g_k)_{k\in \mathbb K}$ such that $g_k$ is in an admissible set which we denote by ${\cal G}_k$. We use $\cal G$ to denote the collection of admissible sets $({\cal G}_k)_{k\in \mathbb K}$

\subsection{Dynamic programming}
Optimal solutions can be obtained through dynamic programming.
 The value of the goal function after time $k$ when starting from state value $x_k$ under policy $g$ is represented by the value functions $V_k^g:{\mathbb R}\to{\mathbb R}$ for all $k\in \mathbb K$ and $g\in\cal G$ 
\begin{equation}
V_k^g(x) = \mathbb{E}\left.\left[ J(x_T^g)-\sum_{j=k}^{T-1} L_j(g_j(x_j^g)) \Delta t   \ \right| \ x_k^g=x\right].     
\end{equation}
The control problem can thus be formulated as finding the control policy resulting in the maximum value
\begin{equation}\label{eq:controlproblem}
V_k^*(x) = \max_{ g_i\in {\cal G}_i\, (k\leq i\leq T-1)} V_k^g(x),
\end{equation} 
for all $k\in \mathbb K$ and for all states $x$ that can be reached on day $k$.

The dynamic programming principle implies that the sequence of functions $V_k^*$ can be found using the backward recursion
\begin{eqnarray}\label{eq:DPg}
    g_k^*(x)&=&\arg\max_{ g_k\in {\cal G}} \mathbb{E}\left[V_{k+1}^*(x_{k+1}^{g_k}) - L_k(g_k(x_k)) \Delta t  \right.  \ | \ x_k=x ],\\
    V_k^*(x)&=& \max_{ g_k\in {\cal G}} \mathbb{E}\left[V_{k+1}^*(x_{k+1}^{g_k}) - L_k(g_k(x_k)) \Delta t  \right.  \ | \  x_k=x ],  \label{eq:DPV}  
\end{eqnarray}
with $V_T^*(x)=J(x)$ and the policy $g^*=\{ g_1^*(x), g_2^*(x),..,g^*_{T-1}(x) \}$ is optimal \cite{Bertsekas2012}.

The conditional expectations depend on the conditional probability of state transitions described by the system model (\ref{eq:model}) \cite{Giuliani2014}, whereas the running costs at time $k$ are known once the policy has been decided. This results in 

\begin{align}\label{eq:expectancyV}
\mathbb{E}[V_{k+1}^*(x_{k+1}^{g_k})\ \mid \ x_k=x] &=\int_{\mathbb R} V^*_{k+1}(y)p_{x_{k+1}^{g_k}|x_k,g_k} (y|x,g_k) dy\\
    \mathbb{E} [L(g(x_k))\ \mid \ x_k=x]&=L(g(x_k)) \label{eq:expectancyL},
\end{align}
and the conditional probabilities of state transitions are described by
\begin{equation}
    p_{x_{k+1}^{g_k}|x_k,g_k} (y|x,g_k )  
    \ = \ \frac{1}{x\sigma   \sqrt{2\pi\Delta t}}\exp{\frac{(y-f_k(x,g_k)\Delta t)^2}{-2\sigma^2 x^2\Delta t}}.
\end{equation}

\subsection{Discretization}
The control problem (\ref{eq:controlproblem}) can be solved by applying the backward induction \eqref{eq:DPg}-\eqref{eq:DPV} on a discretized system. A linear grid was used: $$x(i)=x_{min}+(i-1)\Delta x,\qquad i=1,...,N.$$ 
The probability density of a state jump from $x(i)$ to any $x(j)$, for all $i,j$, is described by a transition matrix $A$, whose elements are defined as
\begin{eqnarray}
      \tilde{A}_{i,j}^g &=& p_{x_{k+1}^{g}|x_k,g_k} (x(j)|x(i),g ) =
    \frac{1}{x(i)\sigma   \sqrt{2\pi\Delta t}}\exp{\frac{(x(j)-f_k(x(i),g)\Delta t)^2}{-2\sigma^2 x(i)^2\Delta t}},\\
   {A}_{i,j}^g &=&\tilde{A}_{i,j}^g /\sum_{j=1}^n\tilde{A}_{i,j}^g.
\end{eqnarray}

We define $V^*_{k,i}=V_k^*(x(i))$ for all $1\leq i\leq N$, which reduces the backward induction in \eqref{eq:DPg}-\eqref{eq:DPV} to
\begin{eqnarray}\label{eq:controlproblemdiscrete}
V^*_{k,i} &=& \max_{ g\in {\cal G}_k}  \sum_{j=1}^N {A}_{i,j}^g (\ V^*_{k+1,j}-L_k(g)\Delta t \ ) ,\\ 
g^*_k &=& \arg\max_{ g\in {\cal G}_k}  \sum_{j=1}^N {A}_{i,j}^g (\ V^*_{k+1,j}-L_k(g)\Delta t \ ).    
\end{eqnarray} 
The values of $V^*_{k,i}$ for different $i$ are collected in a vector $\mathbf{V}_k^*$ and if we build a matrix $\mathbf{A}^g$ from the values of ${A}_{i,j}^g$ the problem can be characterized as
\begin{eqnarray}\label{eq:controlproblemdiscretemat}
V^*_{k} &=& \max_{ g\in {\cal G}_k} \, ( \, \mathbf{A}^g\mathbf{V}_{k+1}^*-L_k(g)\Delta t \, ) ,\\ 
g^*_k &=& \arg\max_{ g\in {\cal G}_k} \, ( \, \mathbf{A}^g\mathbf{V}_{k+1}^*-L_k(g)\Delta t \, )    
\end{eqnarray} 
for $0\leq k\leq T-1$, with $\mathbf{V}_T^*$ the vector with elements $(\mathbf{V}_T^*)_i=V^*(x(i))=J(x(i))$ for all $1\leq i \leq N$.

\subsection{Crop and Greenhouse Model}
We consider the case of lettuce crop (\textit{Lactuca sativa} L.) production inside a greenhouse that is temperature controlled. 
Here $x~[\textrm{kg}\textrm{m}^{-2}]$ is the crop state and the control input at time $k$ (in days)  $u_k=(u^d_k,u^n_k)$ consists of the indoor temperature during the day and during the night, in degrees Celsius. The unit for the scalar plant mass  $x\in \{x(1),x(2),...,x(N)\}$ is kg per square meter and the time step $\Delta t$ represents one day.

\subsubsection{Growth Dynamics}
The crop growth model is based on the single state crop model of \cite{vanHenten1994}, which describes the state dynamics of the dry weight of 16 crops per $\textrm{m}^{2}$ in terms of a function $f_k:\mathbb R^3\to \mathbb R$ that is given by: 
\begin{equation}\label{eq:growthdynamics}
f_k(x_k,u_k^d,u_k^n) =   c_{\rm day} c_\beta(\  c_\alpha F_k^{\rm photo}(x_k,u_k^d,u_k^n,\gamma_k) \ - \  F_k^{\rm resp}(x_k,u_k^d,u_k^n)) \ \ \ [\textrm{kg} \textrm{m}^{-2} \textrm{day}^{-1}] 
\end{equation}
with $x_k$ the crop dry weight on day $k$, $c_{\rm day}=3600\cdot 24$ a scale factor for the number of seconds in a day, and with $c_\beta=0.80$ [-], $c_\alpha=0.68$ [-]. The  photosynthesis as described by the function $F^{\rm photo}:\mathbb R^4\to \mathbb R$ depends, among others, on the sequence $\gamma_k\in {\mathbb R}^{24}$ of average radiation levels $\gamma_{k,h}$ in the Netherlands for each hour $h$ of day $k\in\mathbb K$ which equals
$$
\gamma_{k,h}=c_{\tau,cover}I_{k,h}^{\rm out},
$$
since it is assumed that the measured outdoor radiation $I_{k,h}^{out}$ is partly absorbed or reflected by the cover. Here $c_{\tau,cover}=0.7 [-]$ is the transmissivity of the cover. Empirical observations of  values for $I_{k,h}^{\rm out}$ are available and used as external input. The function $F^{\rm resp}:\mathbb R^3\to \mathbb R$ characterizes respiration. Photosynthesis and respiration are described separately in the following two sub-paragraphs.

\subsubsection{Photosynthesis}
The photosynthesis rate on any given day $k$ is the weighted average over daily and nightly photosynthesis rates which are themselves the average photosynthesis over hours during the day and night respectively\footnote{The criterion to decide whether hour $h$ of day $k$ is day or night is defined as follows: it is night whenever radiation $I_{\rm out}^{k,h}$ during that hour
is lower than the critical value $c_{Ith}=20 \textrm{W}\textrm{m}^{-2}$}. This leads to the specification\footnote{Note that the unit is $s^{-1}$, and that the growth rate is converted to $\rm{day}^{-1}$ by $c_{\rm day}$ in \eqref{eq:growthdynamics}.}
\be
F^{\rm photo}_k(x_k,u_k^d,u_k^n,\gamma_k)&=&  \sfrac{1}{24}\sum_{h=0}^{23}(
 F^{\rm photo}(x_k,T_k^d(u_k^d),\gamma_k^d){\bf 1}_{(k,h){\rm\ is\ day}  }\\&&\qquad
\ + \
 F^{\rm photo}(x_k,T_k^n(u_k^n),\gamma_k^n){\bf 1}_{(k,h){\rm\ is\ night}}) \ [\textrm{kg} \textrm{m}^{-2} \textrm{s}^{-1}].\
\ee
The values are based on the average daytime temperature $T_k(u_k^d)$, the average nighttime temperature $T_k(u_k^n)$ and  average daytime and nighttime radiation levels $\gamma_k^d$ and $\gamma_k^n$ that are defined below; we use these averages instead of hourly values to speed up computations. Temperature averages depend on the control variables for indoor temperatures but also on the  outside temperature for hour $h$ of day $k$, $T_{k,h}^{\rm out}$. It is assumed that on average the indoor temperature is in steady state. Since there is no active cooling (e.g. by pad fan or humidifiers), and since indoor temperature is usually kept equal to, or higher than outside temperature to prevent problems with high relative humidity, it is required that the indoor temperature is equal to, or higher than, the outside temperature. This gives the following averages:
\be
T_k^d(u_k^d) &=& \frac
{\sum_{h=0}^{23} (T_{k,h}^{\rm out}\vee u_k^d){\bf 1}_{(k,h){\rm\, is\ day} }}{\sum_{h=0}^{23} {\bf 1}_{(k,h){\rm\, is\ day}} },\qquad\
\gamma^d_k \ = \ \frac
{\sum_{h=0}^{23} \gamma_{k,h}{\bf 1}_{(k,h){\rm\, is\ day} }}{\sum_{h=0}^{23} {\bf 1}_{(k,h){\rm\, is\ day}} }
\\
T_k^n(u_k^n) &=&\frac
{\sum_{h=0}^{23} (T_{k,h}^{\rm out}\vee u_k^n){\bf 1}_{(k,h){\rm\, is\ night} }}{\sum_{h=0}^{23} {\bf 1}_{(k,h){\rm\, is\ night}} },\qquad
\gamma^n_k \ = \ \frac
{\sum_{h=0}^{23} \gamma_{k,h}{\bf 1}_{(k,h){\rm\, is\ night} }}{\sum_{h=0}^{23} {\bf 1}_{(k,h){\rm\, is\ night}} }.
\ee
Here $a \vee b$ denotes the maximum between $a$ and $b$. The photosynthesis function itself is
\be
F^{\rm photo}(x,T,\gamma) &=&
\phi_{\rm phot,max}(T,\gamma)\cdot(1-e^{-c_k c_{\rm lars}(1-c_{\tau,\rm resp})x}) \ \ [\textrm{kg} \textrm{m}^{-2} \textrm{s}^{-1}]
\ee
with extinction coefficient $c_k=0.90$ [-], shoot leaf area ratio $c_{\rm lars}=62.5$ [m$^2$ kg$^{-1}$], root to total crop dry weight ratio $c_{\tau,\rm resp}=0.07$ [-]. It depends on
\be
\phi_{\rm phot,max}(T,\gamma) &=& \frac{\epsilon(T)c_{\rm par}\gamma \sigma(T)(c_{CO2}-\Gamma(T))}{\epsilon(T)c_{\rm par}\gamma+\sigma(T)(c_{CO2}-\Gamma(T))} \ \ [\textrm{kg} \textrm{m}^{-2} \textrm{s}^{-1}].
\ee
The carbon dioxide compensation point $\Gamma$ in this expression equals 
\be
\Gamma(T) &=& c_{\gamma}\cdot{c_{\rm q10resp}}^{c_{\rm a resp}(T-c_{\rm reftemp})}.
\ee
We have $c_{\gamma}=7.32\cdot10^{-5}$ [kg m$^{-3}$], $c_{\rm q10resp}=2$ [-], $c_{\rm a resp}=0.10$ [$^\textrm{o}\textrm{C}^{-1}$] and $c_{\rm reftemp}=20$ [$^\textrm{o}\textrm{C}$],
with an average density of $\textrm{CO}_\textrm{2}$ in air of $c_{\rm CO2}=0.720\cdot 10^{-3}$ [kg m$^{-3}$] (corresponding to 400 ppm), and
photosynthetically active radiation ratio $c_{\rm par}=0.50$ [-].

The effect of photorespiration on light use efficiency is
\be
\epsilon(T) &=& c_{\epsilon}\frac{c_{CO2}-\Gamma(T)}{ c_{CO2} + 2\Gamma(T) },
\ee
with light use efficiency coefficient $c_{\epsilon} = 17 \cdot 10^{-9}$ [kg J$^{-1}$]. The leaf conductance to carbon dioxide transport function $\sigma$ satisfies
\be
\frac{1}{\sigma(T)} &=&  \frac{1}{c_{\rm bnd}}+\frac{1}{c_{\rm stm}}+\frac{1}{\sigma_{\rm car}(T)}
\ee
with $c_{\rm bnd}=0.004$ [s$^{-1}$] and $c_{\rm stm}=0.007$ [s$^{-1}$]. Further,
\be
\sigma_{\rm car}(T) &=& \ c_{\rm car,2}T^2 + c_{\rm car,1} T + c_{\rm car,0}
\ee
with $c_{car,2}=-5.11\cdot10^{-6} [\textrm{s}^{-2}\, {\,}^\textrm{o}\textrm{C}]$, $c_{car,1}=2.3\cdot10^{-4} [\textrm{s}^{-1}\, {\,}^\textrm{o}\textrm{C}]$ and $c_{car,0}=-6.29\cdot10^{-4} [{\,}^\textrm{o}\textrm{C}]$.

\subsubsection{Respiration Dynamics}
The respiration rate is averaged over the hours in a day, as was the case for photosynthesis,
\be
 F_k^{\rm resp}(x_k,u_k^d,u_k^n) &=&
 \sfrac{1}{24}\sum_{h=0}^{23}(
 F^{\rm resp}(x_k,T_k^d(u_k^d))
{\bf 1}_{(k,h){\rm\ is\ day}  }\\&&\qquad
\ + \
 F^{\rm resp}(x_k,T_k^n(u_k^n)){\bf 1}_{(k,h){\rm\ is\ night}} \ ), \ [\textrm{kg} \textrm{m}^{-2} \textrm{s}^{-1}]
\ee
with
\begin{equation}\label{eq:respiration}
 F^{\rm resp}(x,T) =
x\cdot c_{\rm resp}\ {c_{\rm q10resp}}^{c_{\rm a resp}(T-c_{\rm Tempref})},
\end{equation}
with $c_{\rm Tempref}=25 [{\,}^\textrm{o}\textrm{C}]$, and
$c_{\rm resp}= c_{\rm s resp} (1-c_{\rm \tau resp} )  + c_{\rm r resp}c_{\rm \tau resp},
$ with $c_{\rm \tau resp}=0.07$ [-]. Maintenance respiration for shoot is $c_{\rm s resp}= 3.47\cdot10^{-7}$ [\rm s$^{-1}$], maintenance respiration for root is $ c_{\rm r resp}= 1.16\cdot10^{-7}$ [\rm s$^{-1}$] while as before $c_{\rm q10resp}=2$ [-] and  $c_{\rm a resp}=0.10 [{\,}^\textrm{o}\textrm{C}^{-1}]$.

\subsection{Cost and revenue functions}

\subsubsection{Running costs}
The running costs $L_k$ on day $k$ depend on the amount of heating required. This is a function of the incoming outside radiation $I_{k,h}^{out}$ for hour $h$ of day $k$ and the average daytime and nighttime temperatures based on the control set points $u_k^d$ and $u_k^n$.
Running costs also depend on the outside air temperature $T^{\rm out}_{k,h}$, sky temperature $T^{\rm sky}_{k,h}$ and wind speed $V^{\rm wind}_{k,h}$, in a way that is specified below.

Heating is needed during hour $h$ of day $k$ when $u_k>T^{\rm out}_{k,h}+I_{k,h}/Q_{k,h}$. Furthermore, $I_{k,h}=I_{k,h}^{out}\cdot c_{\tau,cover} \cdot c_{sens}.$ Here $c_{sens}=0.3$ $[-]$ is the fraction of radiation that is converted into sensible heat inside the greenhouse, and $Q_{k,h}$ is the specific heat loss through ventilation and transmission heat per Kelvin temperature difference [$\rm W \textrm{m}^{-2} K^{-1}$]. 

It is assumed that both $c_{\tau,cover}$ and $c_{sens}$ remain constant over time. 
Heating costs are calculated separately for the heating demand in daily and nightly hours since at night it is assumed that the thermal screens are closed, to reduce heat loss. 

The heating cost can then be expressed as:
\be
L_k(u_k^d,u_k^n) &=& c_L \sum_{h=0}^{23}\ [ \ ( \
Q_{k,h}(u_{k}^d)\cdot (u_k^d-T^{\rm out}_{k,h})-I_{k,h}\ )^+ \,
{\bf 1}_{(k,h){\ \rm is\ day}}\\
&&\qquad\ \ +  ( \
Q_{k,h}(u_k^n)\cdot (u_k^n-T^{\rm out}_{k,h})-I_{k,h}\ )^+ \,
{\bf 1}_{(k,h){\ \rm is\ night }}\ ]  ,
\ee
in [\euro\ $\rm{day}^{-1}$],
with $c_L = 3.6\,10^{-6}\, c_{\rm pGJ} / c_{\rm eff}$. The number $3.6\,10^{-6}$ represents the conversion from seconds to hours, and from joule to gigajoule ($[{\rm sh}^{-1} \rm{GJ}\,{\rm (J)}^{-1}]$). Furthermore, $ c_{\rm pGJ} =11$ is the estimated price in \euro\ per GJ gross energy\footnote{ This is based on 
a gross natural gas price of 0.04 \euro/kWh = 11\euro/GJ. }, and $c_{\rm eff}=0.90$ is the water sided efficiency of the heater.

The energy flux $Q_{k,h}$ is calculated as follows:

\be
Q_{k,h}(u) &=&   \rho c_p \cdot\Phi_{k,h} +  c_{\rm As}  \frac{ c_{\rm Sc}{\bf 1}_{(k,h){\ \rm is\ night}} + {\bf 1}_{(k,h){\ \rm is\ day}}}  {1/c_{\alpha} + 1/ (\alpha_{k,h}(u)+{\rm \alpha}_{e,k,h})\ },\\
\alpha_{k,h}(u)&=& 100\wedge  4 c_{\epsilon} \, c_{\sigma}\,
 (\, \sfrac{1}{2}u + \sfrac{1}{2}T^{\rm out}_{k,h} + 273.15\, )^3 \ \left| \frac{u-T^{\rm sky}_{k,h}}{u-T^{\rm out}_{k,h}}\right|,\\
 \Phi_{k,h} &=& (0.2\cdot {\bf 1}_{(k,h){\ \rm is\ night}} + {\bf 1}_{(k,h){\ \rm is\ day}})
\, c_{\rm navg}\, c_{\rm H}.
\ee

Here $\rho c_p$ is specific heat of greenhouse air which is approximately 1206 [$\rm J \, m^{-3} \, K^{-1}$] at 20 $^o$C air temperature, $\Phi_{k,h}$ is the specific ventilation rate $[\rm m\, s^{-1}]$, $c_{\rm As}=1.12$ $[-]$ the specific greenhouse cladding area, $c_{Sc}=0.6$ $[-]$ the energy screen factor, $c_{\rm navg}=0.25/3600$ [$s^{-1}$] the average ventilation rate  when the greenhouse is heated, and $c_{\rm H}=6$ $[\textrm{m}]$ the average greenhouse height. The heat conductivity is a sum of effects due to  radiation and convection so 
$c_{\alpha} = c_{\rm \alpha si} + c_{\rm \alpha ci}$ with  $c_{\rm \alpha si}=4.9$ [$\rm W \, m^{-2} \, K^{-1}$] and $c_{\rm \alpha ci}=2.98$ [$\rm W \, m^{-2} \, K^{-1}$]. The constant
 $c_{\epsilon}=0.90$ $[-]$ is the emission constant for the greenhouse cover material for thermal infrared, and $c_{\sigma}=5.67\cdot 10^{-8}$ $[\rm J \,s^{-1}  \,m^{-2} \, K^{-4}]$ is the Stefan Boltzmann constant. Furthermore,
\be
{\rm \alpha}_{e,k,h} &=& 2.5(V^{\rm wind}_{k,h})^{0.8}{\bf 1}_{V^{\rm wind}_{k,h}>4} + (2.8+1.2V^{\rm wind}_{k,h}){\bf 1}_{V^{\rm wind}_{k,h}\leq 4}.
\ee

Here ${\rm \alpha}_{e,k,h}$ $[-]$ represents the convective heat transfer at the outside of a saw-tooth shaped roof  \cite{bot1983,dezwart1996}.

\subsubsection{Revenues}

The revenues associated with crop yield are expressed by the following goal function

\begin{equation}\label{eq:revenues}
J(x)= \ c_{\rm dryfrac}\, c_{\rm price} x  {\bf 1}_{x \in [x^*-\Delta_H,
x^*+\Delta_H]}.
\end{equation}  
    
\begin{figure}
    \centering
    \includegraphics[width=.7\textwidth]{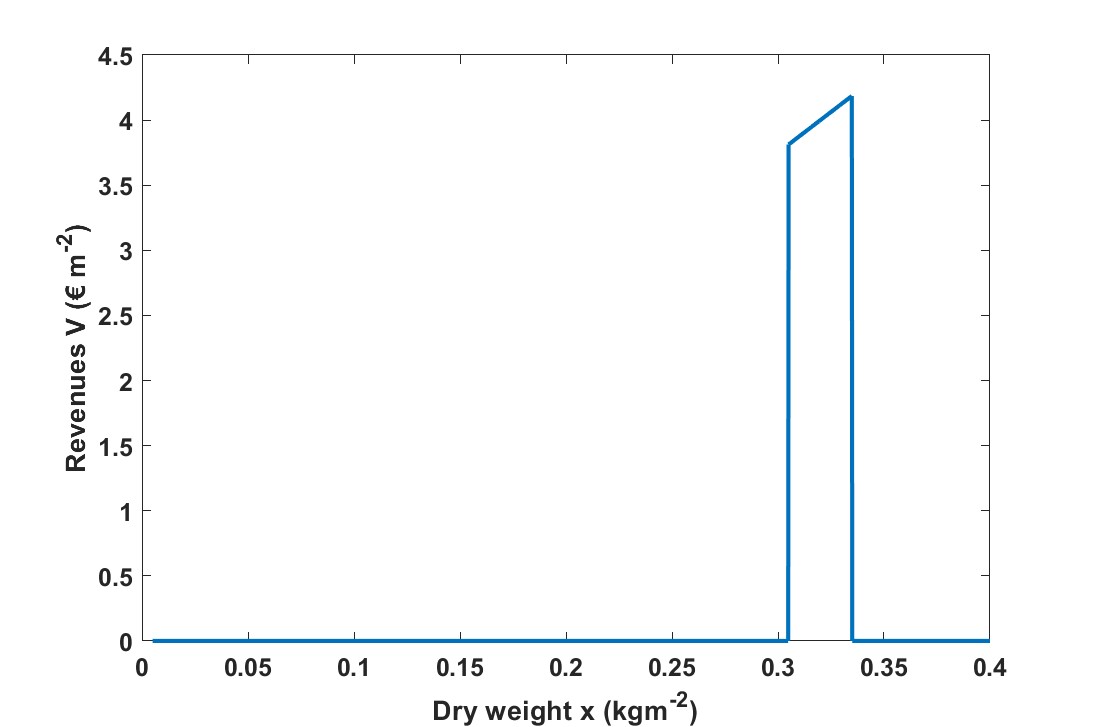}
    \caption{Goal function $J(x)$. The harvested weight is required to fall within a bandwidth of 30 grams around the target weight of 320 grams.}
    \label{fig:goalfunction}
\end{figure}

Here $\Delta_H$ determines the allowed margin of harvested weight, e.g. due to delivery contracts based on consumer demand for product uniformity. The nominal value is $\Delta_H=15 \rm \,g\, m^{-2}$. Further, $c_{\rm dryfrac}=20$ $[-]$ at a dry matter content of 5\% and $c_{\rm price}=1$  [ \euro $\rm \,kg^{-1}]$. The two constants represent the conversion from dry to fresh weight, and the revenues associated with fresh weight, respectively. Assuming a target fresh weight of 400 grams per crop head \cite{KWIN2019}, and 16 plants per $\rm m^{2}$, this comes to a dry weight of $x^*=400/ c_{\rm dryfrac}\cdot 16= 320 \, [\rm g \,m^{-2}]$. The revenue function is displayed in Figure \ref{fig:goalfunction}. 

\subsection{Visualization}

The outcomes of the dynamic programming algorithm to find the optimal policy form a 4-tuple $(u^{d*}_k(x),\, u^{n*}_k(x),\, f^*_k(x),\, V^*_k(x) )$ for every day $k\in\mathbb K$ and every possible value of the state $x$. Here $f^*_k(x)$ is defined as the growth according to the optimal policy given the current state, so $f^*_k(x)=f(x,g^*_k(x))$. 

In the following section we show plots of these four characteristics of the optimal day and night temperatures, optimal plant growth  and optimal revenues to create some intuition for the rationale behind the optimized policies.

\subsection{Computational settings}
\begin{itemize}
    \item All numerical computations were carried out within a Matlab environment.
The computational grid used was $N \times T =2000 \times 40$. That is, the state was discretized into $N=2000$ parts, and the time horizon was $T=40$ days, and $\Delta t =$ 1 day.



\item The state space range is $[5,400]~[\textrm{g} \,\textrm{m}^{-2}]$, and the initial state in the dynamic simulations was $x_0=5\, [\textrm{g}\,\textrm{m}^{-2}]$ unless stated otherwise.
\item The optimization algorithm to solve equation (\ref{eq:controlproblemdiscretemat}) was carried out with the 'fmincon' routine.  
The set of admissible control values is ${\cal G}=[5,10]\times [5,20]$, in $^\textrm{o} \textrm{C}$. This means that the night temperature set point is within the range $[5,10]$ $^\textrm{o} \textrm{C}$, and the day set point is within the range $[5,20]$ $^\textrm{o} \textrm{C}$.
\item The weather data was retrieved for 2014 from Bleiswijk, The Netherlands (52°02N 4°32E), as was done in a previous study \cite{Katzin2020}. 
\item For each selected day, the following period of 40 days is assumed to have the same daily weather pattern as the selected day. This simplification allows us to examine the control strategy based on solely crop state and time before harvest, without changes in weather. The uncertainty as a result of weather forecast errors is assumed to be accounted for by the introduction of uncertainty in the state process. 
\item The maximum initial value was defined as $\mathop{\max}_k V^*_{k,1}$, the maximum expected net revenue at the start of production cycle, with initial value $x=x_0$. 
\end{itemize}

\subsubsection{Dynamic stochastic controller}
The following settings were used to investigate the system dynamics, control policy, and performance. 
\begin{itemize}
\item The nominal weather pattern has relatively much light, and low temperatures (on day 79 in 2014, mean temperature was 4.6 $[^\textrm{o} \textrm{C}]$, and the mean light intensity was $190~[\textrm{W}\textrm{m}^{-2}]$), which creates a heat limitation in the greenhouse. That is, extra heating during the day will result in increased crop growth. 
The control laws and corresponding performance have been computed for two other weather patterns: for day 5 of the year (with mean temperature  -4 $[^\textrm{o} \textrm{C}]$, and mean light intensity $53~[\textrm{W}\textrm{m}^{-2}]$), and for day 187 (with mean temperature  18 $[^\textrm{o} \textrm{C}]$, and mean light intensity $330~[\textrm{W}\textrm{m}^{-2}]$).
\item The nominal level of the parameter $\sigma$ which determines the uncertainty in the state dynamics was $\sigma^2=10^{-4}~[\rm kg^{-1} m^2 day^{-1}]$.

\end{itemize}

\subsubsection{Comparison to other controllers}
The following settings and methods were used to compare the dynamic stochastic controller to two other controllers.
\begin{itemize}
    \item The static stochastic controller input $(u^d,u^n)$, as a function of starting time $t_0$ was determined by carrying out the optimization in  equation (\ref{eq:controlproblem_0}), using the 'fmincon' algorithm in Matlab  for varying starting times, and for fixed initial weight $x_0$.  
    \item A dynamic deterministic optimal controller was found by applying the dynamic stochastic control algorithm at nominal settings, with a variance that was reduced with a factor 50 ($\sigma^2=2 \cdot 10^{-6}$). Strictly speaking this is still a stochastic controller, but it approaches a deterministic one in the sense that it is designed under the assumption that uncertainty is  almost absent (as can be seen in figure \ref{fig:dyn_det}) and therefore deemed negligible.  
    \item The controllers were compared to the dynamic stochastic controller by evaluating the expected net revenues at the optimal starting time. 
    \item A separate sensitivity analysis was carried out for each of four relevant parameters. The performance of the controllers was compared for different levels of uncertainty $\sigma^2$, harvest margin $\Delta H$, starting day, and initial weight $x_0$. The performance was represented by the maximal initial expected revenue $\mathop{\max}_k V^*_{k,1}$ for variations in uncertainty and harvest margin. For variations in starting day and starting weight, the performance was represented by the value ($V^*_{k,i}$ for $k,i$ corresponding to the selected starting day and starting weight $x(i)$). For each simulation the deterministic controller was designed assuming an extremely small value of $\sigma^2$ (as described above), and the performance was calculated using the actual value. For each selected parameter value, and for each controller, the optimal control policy for that specific parameter value was applied.
\end{itemize}

\section{Results}

\begin{figure}
    \centering
    \includegraphics[width=.8\textwidth]{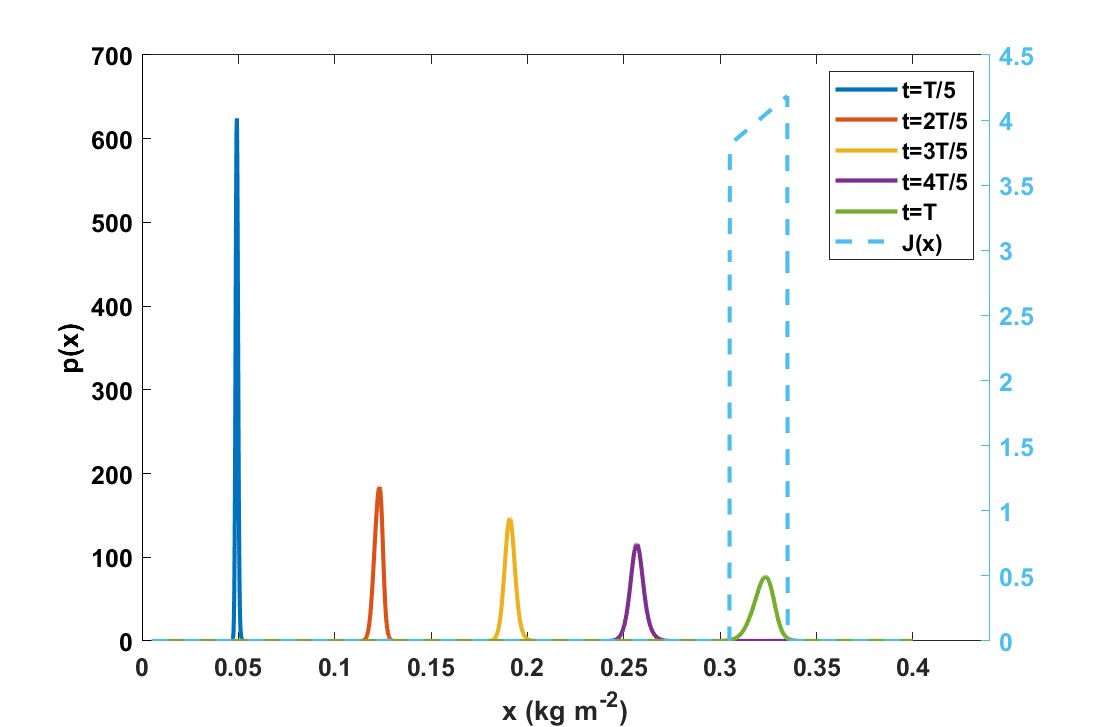}
    \caption{The probability density dynamics of crop weight at 5 equal  time intervals before the final time $T$. The dispersion is small compared to the allowed harvest margin, so in this case there is hardly any risk that the harvest weight will end up outside the allowed weight margin (dotted line).} 
      \label{fig:nominaldynamics}
\end{figure}

The controlled dynamics of the probability density of crop weight for different times after the starting time $t=0$ can be found in Figure \ref{fig:nominaldynamics}, together with the goal function.  
Figure \ref{fig:controllaw} shows 1) the control strategy (daily and nightly temperature set-points), 2) resulting growth rates, and 3) the associated value function $V$ (expected net revenue) for all $x_k$ and $k$ in four heat maps with color scales.

\subsection{Controlled dynamics}\label{sec:systemdynamics}

Figure \ref{fig:nominaldynamics} shows how the probability density of crop weight disperses over time. The final dispersion is relatively small compared to the allowed harvest margin, so there is not much risk that the harvest weight will end up outside the allowed margin (dotted line).

The controller balances optimality with robustness towards the effects of uncertainty. An optimal harvest weight would be 345 $\textrm{g}\,\textrm{m}^{-2}$  since this maximizes $J$ but a distribution centered around 320 $\textrm{g}\,\textrm{m}^{-2}$ since this minimizes the risk that crop weight will fall outside the required weight range. The mean of the state at harvest time $x(T)$ that is realized by the control action is $322~\textrm{g}\,\textrm{m}^{-2}$. This is slightly higher than the target weight of $320 \,\textrm{g}\,\textrm{m}^{-2}$, and the risk of a zero payoff is very small.

\subsection{Control action}

 The optimal indoor night and day temperatures show clearly that the decision to heat or not depends on time and state. This suggests that a controller whose actions do not depend on time or state could  be suboptimal, and this will indeed be confirmed in our later discussion of Figures \ref{fig:perf_static} and \ref{fig:dyn_det}. 

\begin{figure}
    \centering
    \includegraphics[width=1.0\textwidth]{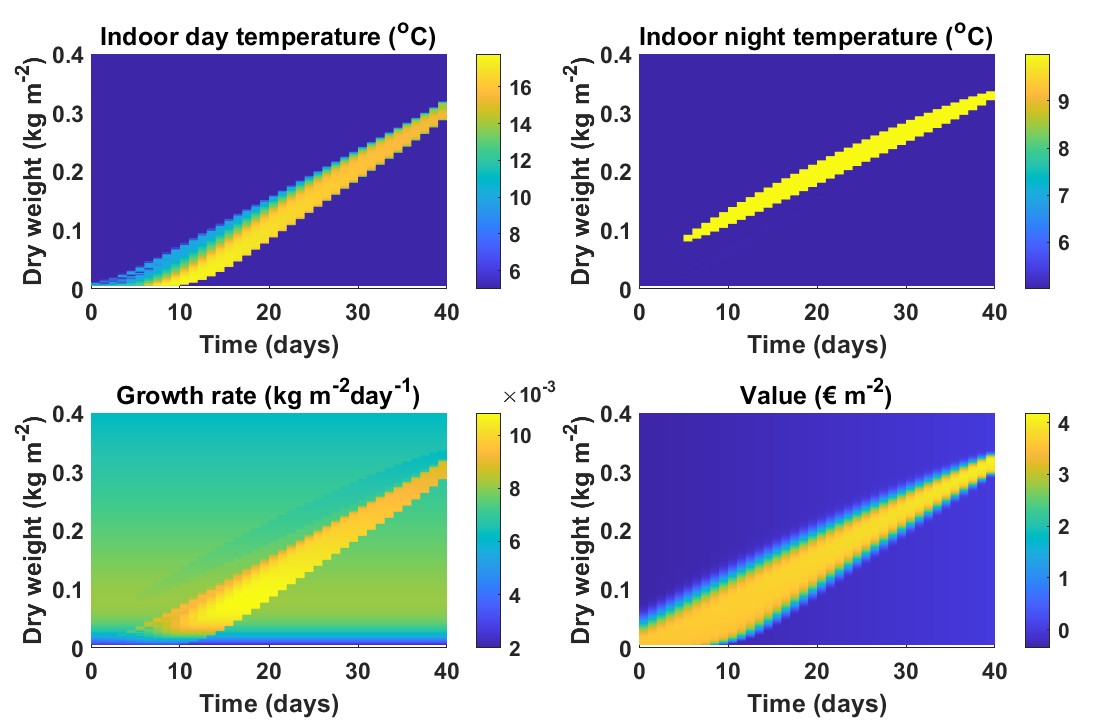}
    \caption{Optimal control policy (top), and corresponding performance (bottom) for day 79 in 2014. Top left: indoor day temperature. Top right: indoor night temperature. The regions in the state space where input is applied form convex bands (in yellow). Bottom left: crop growth rate. Bottom right: value function $V$ (expected net revenues). The controlled regions translate into a single band in the growth rate plot, and also in a band in the value plot. The dispersion of the state, caused by the  state uncertainty, translates into dispersion of the value function.}
    \label{fig:controllaw}
\end{figure}

The regions in the state space with the most control action form convex yellow bands  in figure \ref{fig:controllaw}. In these regions growth rate and consequently expected net revenues are influenced the most. There is little to no control input outside these quite narrow regions, since there is only a single input (heating) that can be exploited. Including more actuators such as lamps and blackout screens might result in larger controlled regions. Energy cost considerations do not seem to influence the sizes of the controlled regions much. In a study not shown here, it was found that decreasing the energy price by a factor 2 barely altered the controlled regions.

\subsection{Controlled growth}

The controlled regions translate into a single band in the growth rate plot where the values are relatively high. For these states it is likely that the crop weight at harvest time will be inside the desired range of harvest weight defined in (\ref{eq:revenues}), i.e. close to the target weight of 320 $\textrm{g} \textrm{m}^{-2}$. Outside the band the value is zero or close to zero, since the harvest weight is expected to fall outside that range. The indoor \textit{day} temperature is used to elevate the growth rate: under high light intensities and low temperatures, growth is heat limited due to inhibition of metabolic processes that convert carbohydrates produced by photosynthesis, to dry matter. The indoor \textit{night} temperatures is used to lower the growth rate (in absence of photosynthesis in the night, a higher temperature will increase the maintenance respiration rate (equation (\ref{eq:respiration})) and thereby decrease growth rate (equation (\ref{eq:growthdynamics})). The regions with elevated day and night temperatures coincide largely with the yellow band in the value plot.

\subsection{Revenues}

The dispersion of the state densities translates into dispersion of the value function, which narrows and brightens as it approaches the time of harvest. An early start of the crop production cycle leads to more dispersion of the bands in the value plot since increasing the length of a production cycle will increase the uncertainty in the outcomes. On average, the value function increases with time, since the remaining cumulative running costs $L$ decrease over time. The values just outside the controlled regions can become slightly negative, due to the constraint that the minimal indoor temperature equals the outdoor temperature. Especially at night, the heat loss due to radiation may require additional heating input, which may result in a negative net revenue. 

The sensitivity of the optimal expected net revenue  with respect to the starting day of the growth cycle is quite limited. This value varies only 1.5 \% over the first 5 days. This shows that the application of a stochastic dynamic controller can compensate for the effects of choosing a different starting day in the production cycle.

\subsection{Influence of weather patterns}

\begin{figure}
    \centering
    \includegraphics[width=0.9\textwidth]{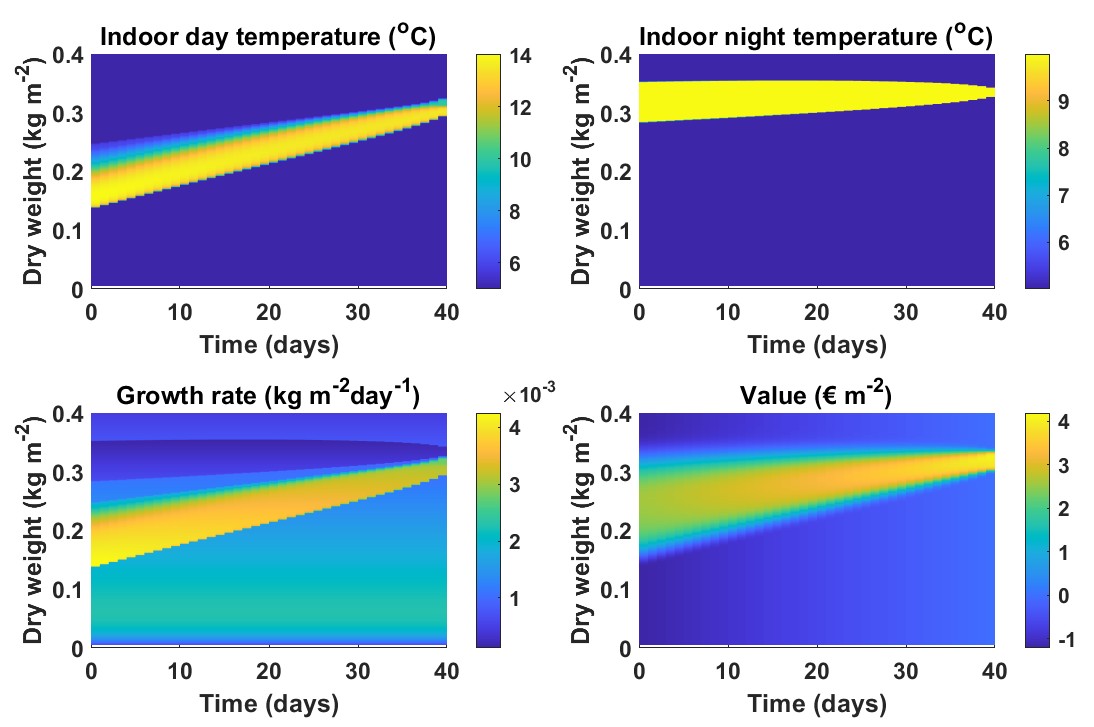}
    
    \includegraphics[width=0.9\textwidth]{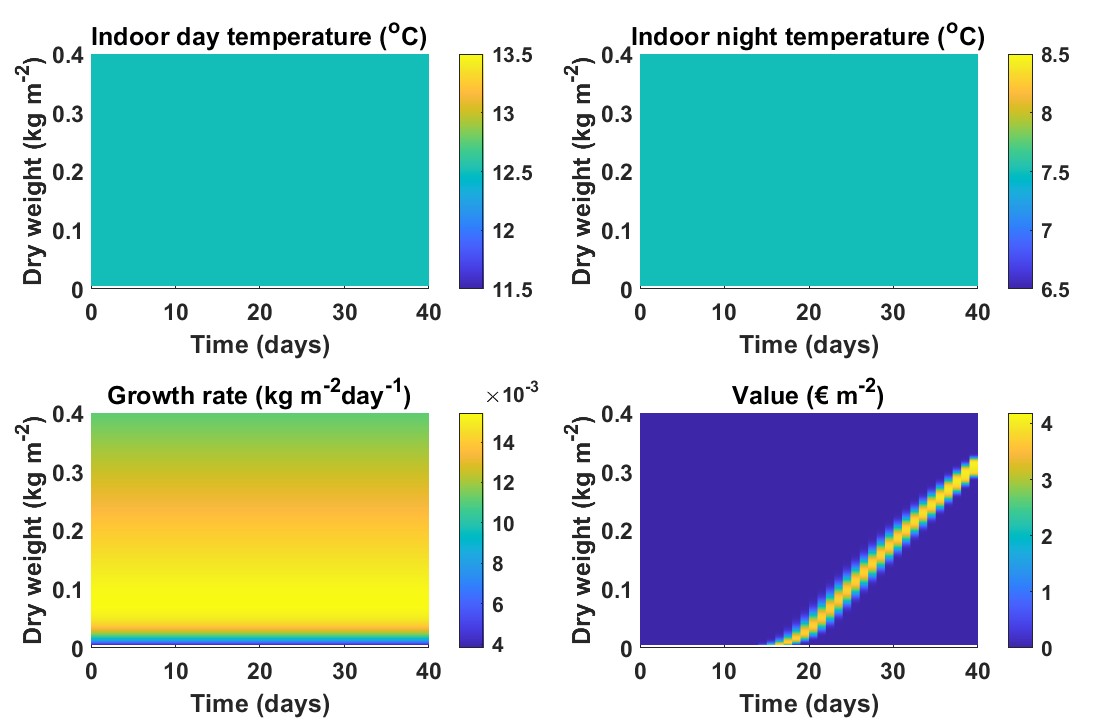}
    \caption{Results for different weather patterns. Top 4 plots: For a weather pattern with low temperature and low light intensity (-4 $^\textrm{o} \textrm{C}$ and $53~\textrm{W}\textrm{m}^{-2}$) the growth rate is more moderate compared to day 79. Bottom 4 plots: For a weather pattern with high temperature and mean light intensity (18 $^\textrm{o} \textrm{C}$ and $330~\textrm{W}\textrm{m}^{-2}$) heating input is not prescribed at any time and for any state. Nonetheless, a high overall growth rate is observed, resulting in a production round that lasts around 24 days, and a steep, lowly dispersed band in the value plot.}
    \label{fig:differentdays}
\end{figure}

 To assess the influence of the weather, optimal control design has been repeated for two other starting days, day 5 and day 187 in 2014. These days were selected since they show weather patterns that differ markedly from that on day 79.  Figure \ref{fig:differentdays}  shows how the control policy depends on the weather.
 
On day 5, there was a low temperature and a low light intensity: the mean temperature was -4 $^\textrm{o} \textrm{C}$ and the mean light intensity was $53~\textrm{W}\textrm{m}^{-2}$.  
The growth rate plot shows that only limited growth is possible. The maximal growth rate is around 4 $\rm g \, m^2 \, day^{-1}$, compared to 10 $\rm g \, m^2 \, day^{-1}$ for day 79, which results in  a more moderate growth curve. The value plot suggests an optimal initial weight of around 0.25 $\textrm{kg}\, \textrm{m}^{-2}$, and during the production cycle of 40 days only 0.07 $\textrm{kg} \, \textrm{m}^{-2}$ of dry weight is gained. However, the costs of starting with such developed plants will be much higher than the costs of a start with smaller plants, and this aspect is not incorporated in the cost function. For a full performance assessment for different weather types, one may need to adapt the duration of the production round to allow a full crop growth cycle for each weather type. Extra costs for running the greenhouse for a longer time, such as staff work hours and rent, would then also have to be incorporated. 

At day 187 growth the mean temperature was 18 $^\textrm{o} \textrm{C}$, and the mean light intensity was $330~\textrm{W}\textrm{m}^{-2}$. This means that growth is more restricted due to the light then as a result of the temperature.  There is therefore no role at all for heating as a means to control growth, and the daily and nightly set points assume the values of the outside temperatures during day and night, respectively. Despite it being  light limited, the maximal growth rate is relatively high: around 14 $\textrm{g}\textrm{m}^{-2} \textrm{day}^{-1}$, compared to 10 $\textrm{g}\textrm{m}^{-2} \textrm{day}^{-1}$ for day 79, due to the high mean light intensity. The value plot indicates an optimal state trajectory with a steep growth curve, and a production round that only lasts around 24 days. The thin yellow band in the value plot can be explained by low dispersion of the state due to a short production round, combined with the absence of feedback control to steer suboptimal outcomes towards better values. 

\section{Comparison with other controllers}\label{sec:othercontrollers}

The stochastic control policy is state and time dependent. Its performance is compared to a static control policy that does not depend on state or time, but is optimized for stochastic state dynamics. The second comparison is against a dynamic deterministic controller, that is optimized with respect to state and time, but does not deal with uncertainty in the state.  

 \subsection{Comparison of controlled state dynamics}

\begin{figure}
    \centering
    \includegraphics[width=0.48\textwidth]{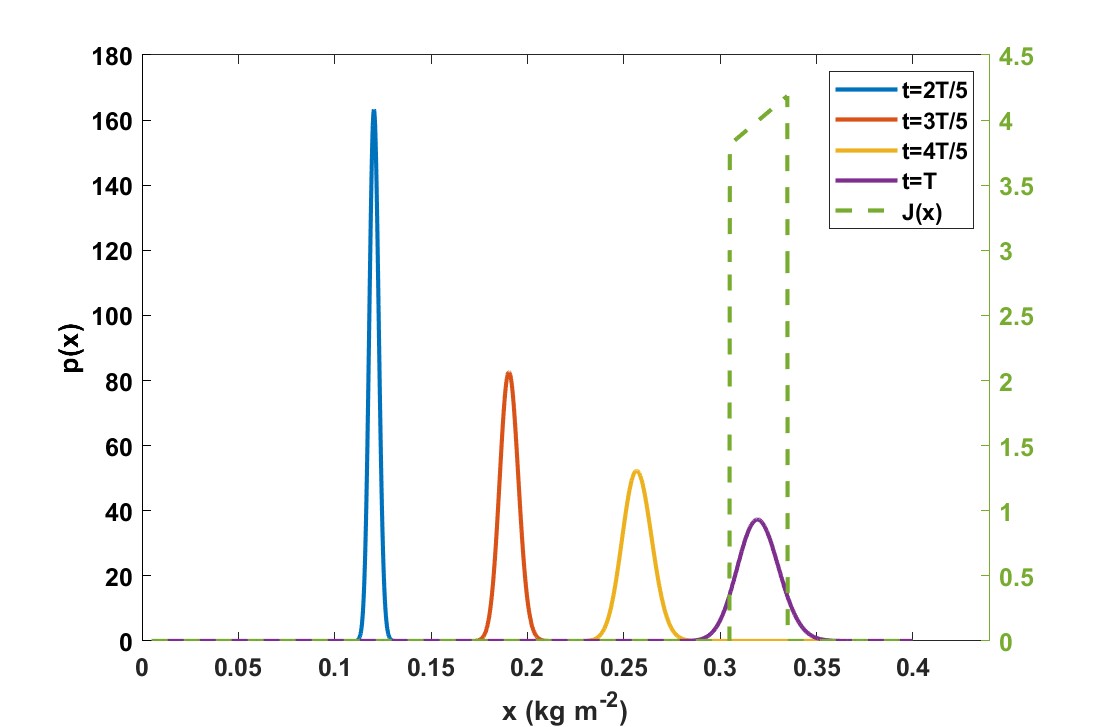}
        \includegraphics[width=0.48\textwidth]{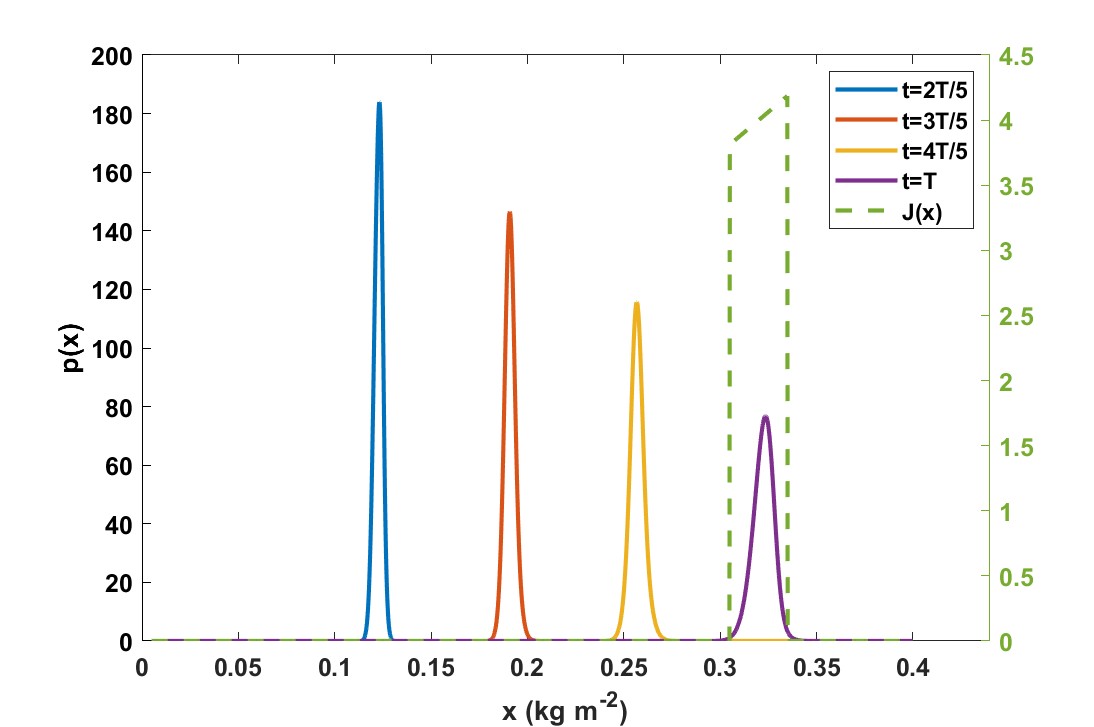}
    \caption{Dynamics of probability through static stochastic control (left), and dynamic stochastic control (right). With static stochastic control the standard deviation at harvest time is much larger compared to dynamic stochastic control.}
    \label{fig:perf_static}
\end{figure}
Figure \ref{fig:perf_static} shows The dynamics of state probabilities with a static stochastic controller, compared to the dynamic stochastic controller for the nominal case, starting from day 1. 
It can be seen that the uncertainty  is considerably larger when this static control is used.  At harvest time the standard deviation of the state is almost twice as large compared to the dynamic stochastic controller ($11$ versus $5.5~ \textrm{g} \textrm{m}^{-2}$), because the latter actively steers the weight towards the optimal trajectory by adjusting the input whereas the static controller does not have such a compensation mechanism. Furthermore, static control results in a smaller maximal value at the starting time compared to the dynamic stochastic controller ($3.0$ versus $3.6$  \euro $\,\textrm{m}^{-2}$ ).

\begin{figure}
    \centering
    \includegraphics[width=0.48\textwidth]{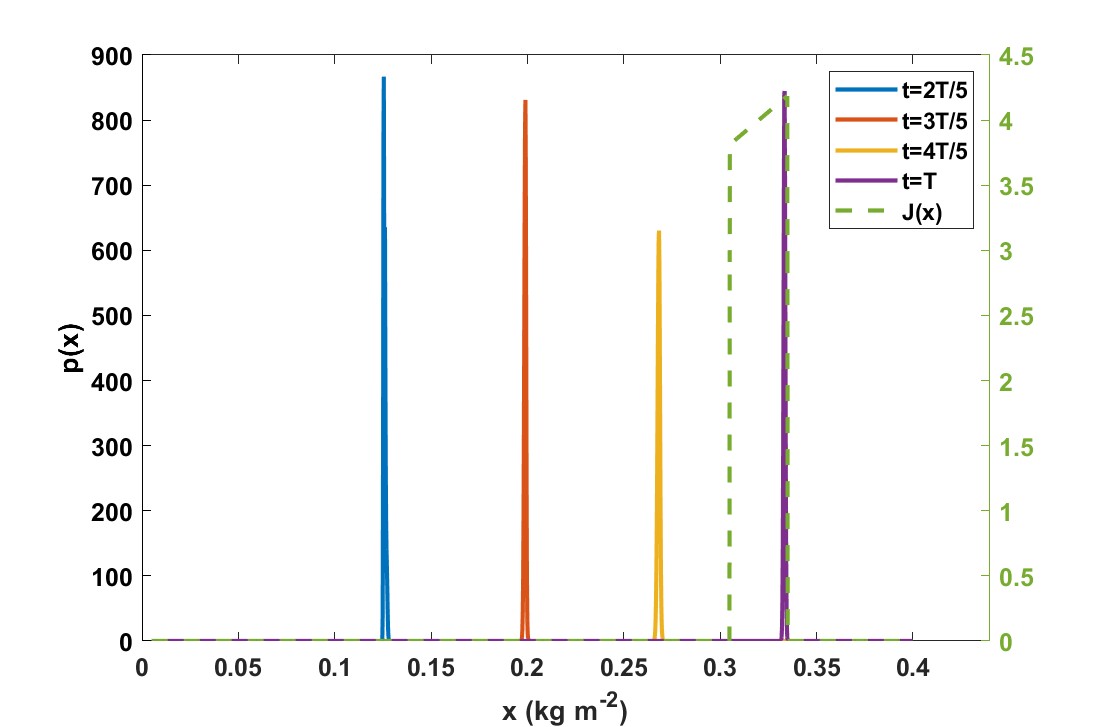}
        \includegraphics[width=0.48\textwidth]{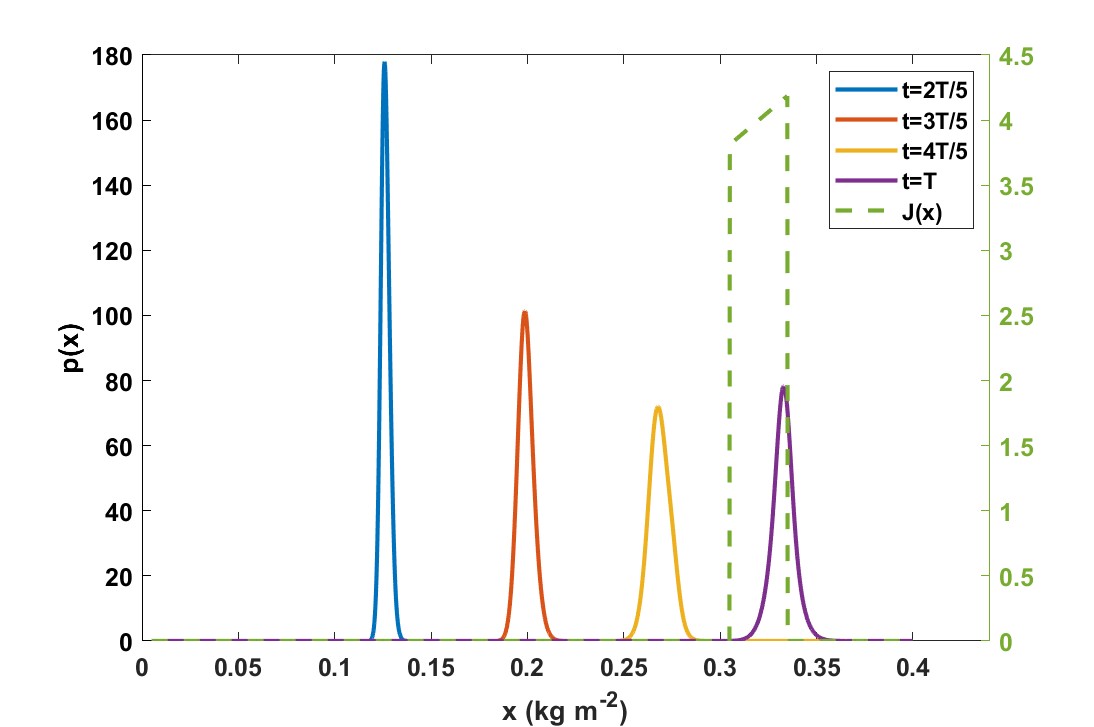} 
    \caption{Dynamics of probability using deterministic control. Left: deterministic controller applied to deterministic system. Right: dynamic deterministic controller applied to stochastic system.}
    \label{fig:dyn_det}
\end{figure}

 \subsection{Comparison of controller action}
 
Figure \ref{fig:dyn_det} shows the state dynamics when a  dynamic but deterministic controller is used.  
When the controller is applied to the (almost) deterministic system that it is designed for, it  steers the state precisely towards the value which maximizes revenues. However, when it is applied to the stochastic system, the standard deviation at harvest time is larger than the one generated by the dynamic stochastic controller: $6.0$ instead of $5.5~ \textrm{g} \textrm{m}^{-2}$. The state that maximizes $J$ is theoretically optimal but not very robust, because even a small positive deviation from this optimum will result in zero revenues from harvesting. Since the controller is not designed to deal with uncertainty, it does not balance optimality in expectation with robustness in the face of uncertainty.

\begin{figure}[!ht]
    \centering
    \includegraphics[width=0.99 \textwidth]{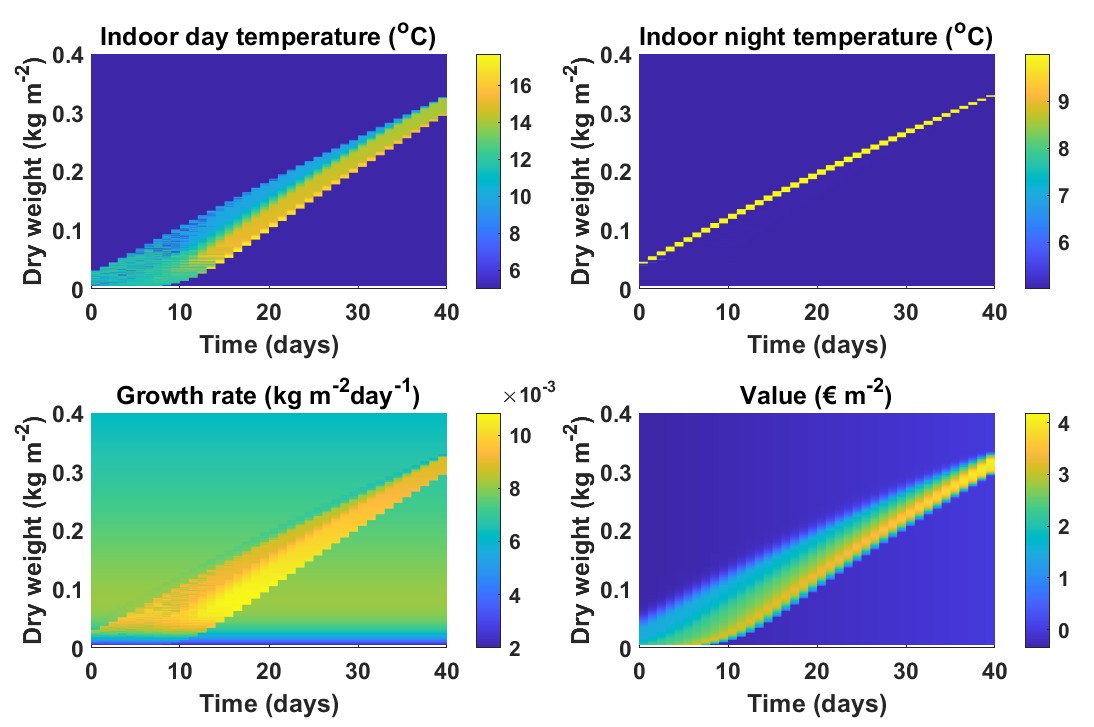}   
    \caption{ Optimal policy, and corresponding performance for the dynamic deterministic controller. Comparison with the dynamic stochastic controller shown in a previous figure shows that the controlled region is more time and state specific. The value plot shows a relatively small band with high values, indicating a limited robustness against deviations from the optimal state trajectory.}
    \label{fig:perf_deterministic}
\end{figure}

Figure \ref{fig:perf_deterministic} shows the control policy and performance of the deterministic controller.
Comparison with Figure \ref{fig:controllaw} shows that for indoor day temperature, the size of the controlled region is now a bit larger.  The magnitude of the control input is a subtle function of time and state, as can be seen by the nuanced color changes within the control band. This may be caused by the assumption that the dynamics are (almost) precisely known, which allows for more subtle input control.

Applying the deterministic controller to the stochastic system results in a considerable loss of initial value compared to the dynamic stochastic controller ($3.1$ versus $3.6$ \euro $\,\textrm{m}^{-2}$). There is no sharp contrast anymore between regions with strong or weak control action:  the values decline more gradually. The fact that the bandwidth with active control is much smaller than for the dynamic stochastic controller indicates a loss of robustness towards deviations of the state from its optimal trajectory. 

\subsection{Comparison of performance}

\begin{figure}[!ht]
    \centering
    \includegraphics[width=0.99 \textwidth]{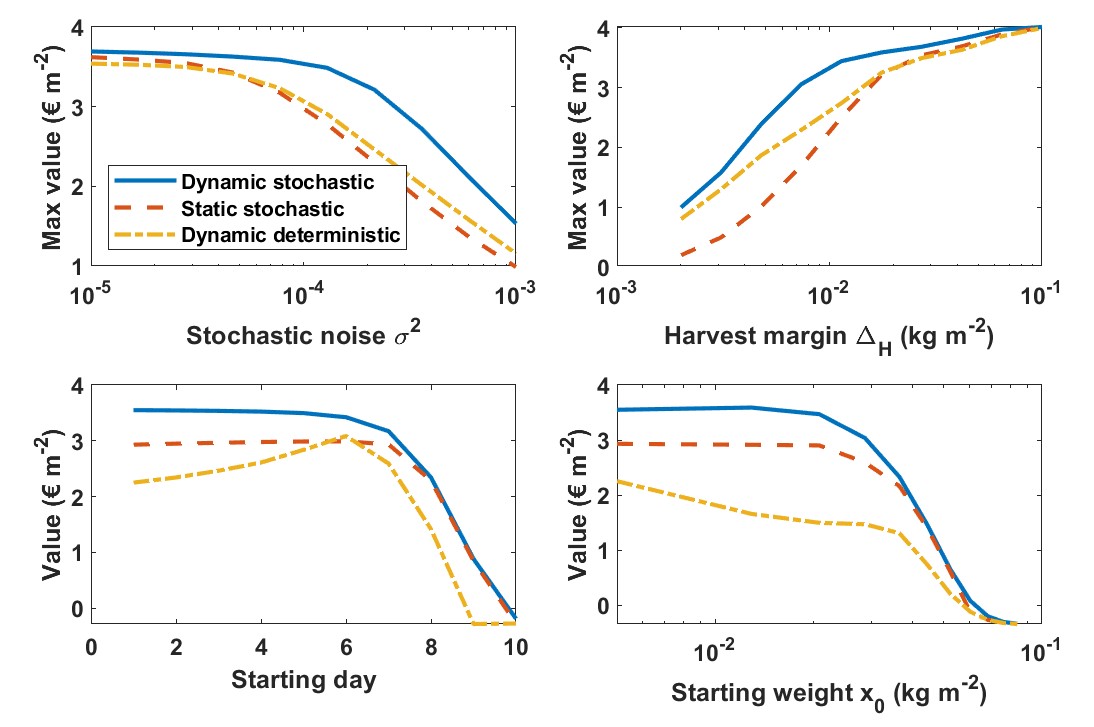}   
    \caption{Sensitivity of performance for a change in $\sigma^2$, harvest margin $\Delta H$, starting day $t_0$, and starting weight $x_0$. For all choices of these parameters, the dynamic stochastic controller outperforms the static stochastic and the dynamic deterministic controllers. For configurations that allow high performance (very little uncertainty, large harvest margin) the differences are relatively small.}
    \label{fig:sensitivity}
\end{figure}

Figure \ref{fig:sensitivity} shows the sensitivity of the performance, in terms of  net revenue, for changes in the following configuration parameters: the level of uncertainty  $\sigma^2$, the harvest margin $\Delta H$, the starting day $t_0$ and the starting weight $x_0$.

The dynamic stochastic controller outperforms the static stochastic and dynamic deterministic controllers for all parameter choices. This was to be expected since for the latter two state dynamics or uncertainty are discarded in the design process.
For some cases, the added value of including dynamics or uncertainty turns out to be substantial, but the increase in value varies with the value of the configuration parameters.  

The upper left plot reveals that  the dynamic stochastic controller outperforms the other controllers by a large margin for a substantial range of levels of uncertainty. If uncertainty is small we observe similar performance under dynamic deterministic and dynamic stochastic control laws. For very small values of $\sigma^2$ we also observe limited performance differences between the static stochastic and dynamic stochastic controllers, indicating that in this case the state dynamics do not play a large role. 
When uncertainty is large, the performance deteriorates for all controllers, since the dispersion in the state density causes the standard deviation at harvest time to be  larger than the allowed harvest margin.

The upper right plot shows performance as a function of allowed harvest margins. For values above $300~ \textrm{g} \textrm{m}^{-2}$ all controllers show more or less the same performance. As explained in \ref{sec:othercontrollers},  dynamic feedback and risk mitigation lead to the strongest improvements under more precise harvest requirements. A decrease in harvest margin to levels below $100~ \textrm{g} \textrm{m}^{-2}$ results in a more or less linear decrease in performance for all controllers. 

The lower left plot shows the sensitivity with respect to a change in the start of the production round. It shows that the deterministic controller is not very robust in the face of such changes. This confirms the earlier observation  from Figure \ref{fig:perf_deterministic} that the bandwidth of high revenue values is very small.

The results in the lower right plot also show the severe limitations of the dynamic deterministic controller for different starting weights: losses may be as high as $50\%$ for some cases. This confirms what can be seen from the values along the vertical axis at the starting day  in Figure \ref{fig:perf_deterministic}.

\section{Conclusions and suggestions for further research}

The dynamic stochastic controller that we designed for a greenhouse model with uncertain dynamics outperformed the stochastic but static controller and the dynamic but deterministic controller considerably in terms of maximal expected net revenue (by 19 \% and 15 \% respectively). This testifies to the value of incorporating knowledge of time and state variables in greenhouse control systems, to allow time-varying feedback and the mitigation of uncertainty. 

The dynamic stochastic controller also achieved better state precision at harvest time (50 \% and 8\% less standard deviation, respectively) than the static controller. This suggests that time-varying feedback may increase state precision dramatically, and that taking uncertainty into account may increase state precision considerably. 

The sensitivity analysis results indicate that the added value is largest in the case of large uncertainty in state dynamics, high required precision in harvest weight, and possibly large deviations from the optimal starting day and the optimal starting weight. 
As discussed in the introduction, these cases are of particular relevance in the modern greenhouse industry, where prediction uncertainty is generally large. 


\subsection{Crop physiology}

By visually linking together the optimal control policy, subsequent crop dynamics and performance in terms of expected net revenues as a function of time and state (figure \ref{fig:controllaw}), transparency is created in the sense that some mechanisms that underlie the outcomes of a rather complex optimization algorithm become intuitively clear. In particular, the comparison between control input, crop dynamics, and the value function reveals how temperature adjustments are used to steer the weight in case it deviates from its optimal trajectory towards a desirable harvest weight. 

It is particularly interesting to see that an elevated indoor temperature is prescribed to accelerate the crop growth rate during daytime and to decelerate it during nighttime. Heating during the day to accelerate the growth makes sense physiologically, since increasing assimilate production and conversion to cell matter requires heat and light. Heating during the night to decelerate the growth makes sense as well, since without light there will be no dry matter production, so assimilate loss will increase due to increased maintenance respiration. 

However, it is well known that an imbalance between heat and light input can result in thick yellow leaves due to starch buildup in case of heat limitation, or create long thin leaves in case of light limitation. One may want to avoid  violations to the requirements regarding the heat-light balance  through refinements in the definition of the control input domain $\cal G$. 
This may also help to define further model improvements in the form of model extensions, such as added state variables or control variables. 

\subsection{Computational demand}

Model extensions come at a computational price as a result of the so-called curse of dimensionality. A possible way to mitigate higher computational costs is to employ a reinforcement learning strategy \cite{Tchamitchian2005, zhang2021, Ajagekar2022}, where an optimal policy is constructed by exploring only a fraction of all 
possible states. As can be see in the control map (e.g. figure \ref{fig:controllaw}), a large fraction of the possible states is unlikely to be realized (large weights near the start, and small weights near the end of the cultivation round), hence an optimal policy does not necessarily need to cover those instances. 


\subsection{Refinement of the model for uncertainty}

In our greenhouse specification, all uncertainty sources are indirectly modelled via the stochastic dynamics of the 
state. In order to disentangle the roles of the various sources of uncertainty, a promising approach could be to extend the framework with a state estimator such as the Kalman filter, to incorporate the fact that in practice it may be impossible to observe the state exactly.
In this paper, the control policy is designed based on the assumption that the crop state can be measured online and without error. In practice, online crop weight measurement is still not standard, even in sophisticated greenhouse systems, since it requires costly weight sensors or the employment of state of the art computer vision via cameras to estimate weight in a smart way. Promising results have been achieved e.g. via point cloud estimation \cite{Mortensen2018segmentation}, or a neural network approach \cite{zhang2020}. In these studies the measurement errors were quite small, resulting in high $R^2$ values (between $0.8$ and $0.95$) but they may still negatively effect control performance. Hence, it would be useful to investigate the effect of measurement accuracy on control performance in further research.

Another possible extension of the framework would be to incorporate uncertainty in the running costs. Energy prices vary over time in an unpredictable manner and one might want to incorporate the stochastic nature of energy costs in a more realistic model. In a similar vein, modelling uncertainty in the weather may help to design controllers which are even more robust to the changing circumstances in which they have to operate. 

\medskip

We believe that these and the other suggestions for further studies may offer valuable additional insights, but we expect that they will leave the conclusions of the present study intact:  risk mitigation and a time-varying feedback control policy can lead to considerable improvements for greenhouse management in an uncertain environment.

\section*{Acknowledgements}
We would like to thank Christian Hamster (Wageningen University) for useful discussions. 



\bibliography{bibfileStoch}

\begin{thebibliography}{10}

\bibitem{Poorter2013}
Hendrik Poorter, Niels~PR Anten, and Leo~FM Marcelis.
\newblock Physiological mechanisms in plant growth models: do we need a
  supra-cellular systems biology approach?
\newblock {\em Plant, cell \& environment}, 36(9):1673--1690, 2013.

\bibitem{Righini2020}
Isabella Righini, Bram Vanthoor, Mich{\`e}l~J Verheul, Muhammad Naseer, Henk
  Maessen, Tomas Persson, and Cecilia Stanghellini.
\newblock A greenhouse climate-yield model focussing on additional light, heat
  harvesting and its validation.
\newblock {\em Biosystems Engineering}, 194:1--15, 2020.

\bibitem{balendonck2010}
Jos Balendonck, EA~Van~Os, Rob van~der Schoor, BAJ Van~Tuijl, and LCP Keizer.
\newblock Monitoring spatial and temporal distribution of temperature and
  relative humidity in greenhouses based on wireless sensor technology.
\newblock In {\em International Conference on Agricultural Engineering-AgEng},
  pages 443--452, 2010.

\bibitem{Speetjens2009}
SL~Speetjens, JD~Stigter, and G~Van~Straten.
\newblock Towards an adaptive model for greenhouse control.
\newblock {\em Computers and electronics in agriculture}, 67(1-2):1--8, 2009.

\bibitem{Hameed2010}
Ibrahim~A Hameed.
\newblock Using the extended kalman filter to improve the efficiency of
  greenhouse climate control.
\newblock {\em International Journal of innovative Computing, information and
  control}, 6(6):2671--2680, 2010.

\bibitem{vanMourik2019}
S.~van Mourik, van Beveren P. J.~M., I.~L. López-Cruz, and E.~J. van Henten.
\newblock Improving climate monitoring in greenhouse cultivation via model
  based filtering.
\newblock {\em Biosystems Engineering, vol. 181, pp. 40-51}, 2019.

\bibitem{Kuijpers2021}
W.~J.~P. Kuijpers, D.~Katzin, S.~van Mourik, D.~J. Antunes, S.~Hemming, and
  M.~J.~G. van~de Molengraft.
\newblock Lighting systems and strategies compared in an optimally controlled
  greenhouse.
\newblock {\em Biosystems Engineering, vol. 202, pp. 195-216}, 2021.

\bibitem{Su2021}
Yuanping Su, Lihong Xu, and Erik~D. Goodman.
\newblock Multi-layer hierarchical optimisation of greenhouse climate setpoints
  for energy conservation and improvement of crop yield.
\newblock {\em Biosystems Engineering}, 205:212--233, 2021.

\bibitem{Payne2021}
Henry~J. Payne, Silke Hemming, Bram~A.P. {van Rens}, Eldert~J. {van Henten},
  and Simon {van Mourik}.
\newblock Quantifying the role of weather forecast error on the uncertainty of
  greenhouse energy prediction and power market trading.
\newblock {\em Biosystems Engineering}, 224:1--15, 2022.

\bibitem{Vazquez2014}
MA~Vazquez-Cruz, R~Guzman-Cruz, IL~Lopez-Cruz, O~Cornejo-Perez,
  I~Torres-Pacheco, and RG~Guevara-Gonzalez.
\newblock Global sensitivity analysis by means of efast and sobol methods and
  calibration of reduced state-variable tomgro model using genetic algorithms.
\newblock {\em Computers and Electronics in Agriculture}, 100:1--12, 2014.

\bibitem{Lopez2018uncertainty}
IL~L{\'o}pez-Cruz, A~Mart{\'\i}nez-Ruiz, A~Ruiz-Garc{\'\i}a, and M~Gallardo.
\newblock Uncertainty analyses of the vegsyst model applied to greenhouse
  crops.
\newblock In {\em XXX International Horticultural Congress IHC2018: III
  International Symposium on Innovation and New Technologies in Protected
  1271}, pages 199--206, 2018.

\bibitem{Guzman2009}
Rosario Guzm{\'a}n-Cruz, R~Casta{\~n}eda-Miranda, Juan~Jos{\'e}
  Garc{\'\i}a-Escalante, IL~L{\'o}pez-Cruz, Alfredo Lara-Herrera, and JI~De~la
  Rosa.
\newblock Calibration of a greenhouse climate model using evolutionary
  algorithms.
\newblock {\em Biosystems engineering}, 104(1):135--142, 2009.

\bibitem{Oldewurtel2013}
Frauke Oldewurtel, Colin~Neil Jones, Alessandra Parisio, and Manfred Morari.
\newblock Stochastic model predictive control for building climate control.
\newblock {\em IEEE Transactions on Control Systems Technology},
  22(3):1198--1205, 2013.

\bibitem{boucherie2017}
Richard~J Boucherie and Nico~M Van~Dijk.
\newblock {\em Markov decision processes in practice}.
\newblock Springer, 2017.

\bibitem{Ross2014}
Sheldon~M Ross.
\newblock {\em Introduction to stochastic dynamic programming}.
\newblock Academic press, 2014.

\bibitem{Huong2018}
Truong~Thu Huong, Nguyen~Huu Thanh, Nguyen~Thi Van, Nguyen~Tien Dat, Nguyen
  Van~Long, and Alan Marshall.
\newblock Water and energy-efficient irrigation based on markov decision model
  for precision agriculture.
\newblock In {\em 2018 IEEE Seventh International Conference on Communications
  and Electronics (ICCE)}, pages 51--56. IEEE, 2018.

\bibitem{kovacs2012}
S~Kov{\'a}cs, P~Balogh, et~al.
\newblock Application of the hierarchic markovian decision processes in the
  decision making processes of pig keeping.
\newblock {\em Agr{\'a}rinformatika Foly{\'o}irat}, 3(2):37--49, 2012.

\bibitem{Onstad1985}
David~W Onstad and Rudy Rabbinge.
\newblock Dynamic programming and the computation of economic injury levels for
  crop disease control.
\newblock {\em Agricultural Systems}, 18(4):207--226, 1985.

\bibitem{sells1995}
JE~Sells.
\newblock Optimising weed management using stochastic dynamic programming to
  take account of uncertain herbicide performance.
\newblock {\em Agricultural Systems}, 48(3):271--296, 1995.

\bibitem{Boussios2019}
David Boussios, Paul~V Preckel, Yigezu~A Yigezu, Prakash~N Dixit, Samia
  Akroush, Hatem~Cheikh M'hamed, Mohamed Annabi, Aden Aw-Hassan, Yahya
  Shakatreh, Omar Abdel~Hadi, et~al.
\newblock Modeling producer responses with dynamic programming: A case for
  adaptive crop management.
\newblock {\em Agricultural Economics}, 50(1):101--111, 2019.

\bibitem{Kuijpers2022}
Wouter~J.P. Kuijpers, Duarte~J. Antunes, Simon {van Mourik}, Eldert~J. {van
  Henten}, and Marinus~J.G. {van de Molengraft}.
\newblock Weather forecast error modelling and performance analysis of
  automatic greenhouse climate control.
\newblock {\em Biosystems Engineering}, 214:207--229, 2022.

\bibitem{Della2019}
A~Della~Noce, M~Carrier, and P-H Courn{\`e}de.
\newblock Optimal control of non-smooth greenhouse models.
\newblock In {\em International Symposium on Advanced Technologies and
  Management for Innovative Greenhouses: GreenSys2019 1296}, pages 125--132,
  2019.

\bibitem{Zhuang2018}
Peng Zhuang, Hao Liang, and Mitchell Pomphrey.
\newblock Stochastic multi-timescale energy management of greenhouses with
  renewable energy sources.
\newblock {\em IEEE Transactions on Sustainable Energy}, 10(2):905--917, 2018.

\bibitem{garcia2023multi}
Francisco Garc{\'\i}a-Ma{\~n}as, Francisco Rodr{\'\i}guez, Manuel Berenguel,
  and Jos{\'e}~Mar{\'\i}a Maestre.
\newblock Multi-scenario model predictive control for greenhouse crop
  production considering market price uncertainty.
\newblock {\em IEEE Transactions on Automation Science and Engineering}, 2023.

\bibitem{Kocian2020}
Alexander Kocian, Daniele Massa, Samantha Cannazzaro, Luca Incrocci, Sara
  Di~Lonardo, Paolo Milazzo, and Stefano Chessa.
\newblock Dynamic bayesian network for crop growth prediction in greenhouses.
\newblock {\em Computers and Electronics in Agriculture}, 169:105167, 2020.

\bibitem{Tchamitchian2005}
Marc Tchamitchian, Constantin Kittas, Thomas Bartzanas, and Christos Lykas.
\newblock Daily temperature optimisation in greenhouse by reinforcement
  learning.
\newblock {\em IFAC Proceedings Volumes}, 38(1):131--136, 2005.
\newblock 16th IFAC World Congress.

\bibitem{zhang2021}
Wanpeng Zhang, Xiaoyan Cao, Yao Yao, Zhicheng An, Xi~Xiao, and Dijun Luo.
\newblock Robust model-based reinforcement learning for autonomous greenhouse
  control.
\newblock In Vineeth~N. Balasubramanian and Ivor Tsang, editors, {\em
  Proceedings of The 13th Asian Conference on Machine Learning}, volume 157 of
  {\em Proceedings of Machine Learning Research}, pages 1208--1223. PMLR,
  17--19 Nov 2021.

\bibitem{Ajagekar2022}
Akshay Ajagekar and Fengqi You.
\newblock Deep reinforcement learning based automatic control in semi-closed
  greenhouse systems.
\newblock {\em IFAC-PapersOnLine}, 55(7):406--411, 2022.
\newblock 13th IFAC Symposium on Dynamics and Control of Process Systems,
  including Biosystems DYCOPS 2022.

\bibitem{cao2022igrow}
Xiaoyan Cao, Yao Yao, Lanqing Li, Wanpeng Zhang, Zhicheng An, Zhong Zhang,
  Li~Xiao, Shihui Guo, Xiaoyu Cao, and Meihong Wu.
\newblock igrow: A smart agriculture solution to autonomous greenhouse control.
\newblock In {\em Proceedings of the AAAI Conference on Artificial
  Intelligence}, pages 11837--11845, 2022.

\bibitem{Bertsekas2012}
Dimitri Bertsekas.
\newblock {\em Dynamic programming and optimal control: Volume I}, volume~1.
\newblock Athena scientific, 2012.

\bibitem{Giuliani2014}
Matteo Giuliani, Stefano Galelli, and Rodolfo Soncini-Sessa.
\newblock A dimensionality reduction approach for many-objective markov
  decision processes: Application to a water reservoir operation problem.
\newblock {\em Environmental Modelling \& Software}, 57:101--114, 2014.

\bibitem{vanHenten1994}
E.~J. van Henten.
\newblock {\em Greenhouse climate management: an optimal control approach}.
\newblock PhD thesis, University Wageningen, 1994.

\bibitem{bot1983}
Gerardus~PA Bot.
\newblock {\em Greenhouse climate: from physical processes to a dynamic model}.
\newblock Wageningen University and Research, 1983.

\bibitem{dezwart1996}
HF~De~Zwart.
\newblock {\em Analyzing energy-saving options in greenhouse cultivation using
  a simulation model}.
\newblock Wageningen University and Research, 1996.

\bibitem{KWIN2019}
M~Raaphorst, J~Benninga, and BA~Eveleens.
\newblock Quantitative information on dutc h greenhouse horticulture 2019.
\newblock {\em Report WPR-898. Wageningen University and Research, The
  Netherlands}, 2019.

\bibitem{Katzin2020}
D.~Katzin, S.~van Mourik, F.~L.~K. Kempkes, and E.~J. van Henten.
\newblock Energy saving measures in optimally controlled greenhouse lettuce
  cultivation.
\newblock {\em International Symposium on Innovation and New Technologies in
  Protected Cultivation}, 2020.

\bibitem{Mortensen2018segmentation}
Anders~Krogh Mortensen, Asher Bender, Brett Whelan, Margaret~M Barbour, Salah
  Sukkarieh, Henrik Karstoft, and Ren{\'e} Gislum.
\newblock Segmentation of lettuce in coloured 3d point clouds for fresh weight
  estimation.
\newblock {\em Computers and Electronics in Agriculture}, 154:373--381, 2018.

\bibitem{zhang2020}
Lingxian Zhang, Zanyu Xu, Dan Xu, Juncheng Ma, Yingyi Chen, and Zetian Fu.
\newblock Growth monitoring of greenhouse lettuce based on a convolutional
  neural network.
\newblock {\em Horticulture research}, 7, 2020.

\end{thebibliography}

\end{document}